
\tolerance=10000
\raggedbottom

\baselineskip=15pt
\parskip=1\jot

\def\sk{\vskip 3\jot}

\def\heading#1{\vskip3\jot{\noindent\bf #1}}
\def\label#1{{\noindent\it #1}}
\def\QED{\hbox{\rlap{$\sqcap$}$\sqcup$}}


\def\ref#1;#2;#3;#4;#5.{\item{[#1]} #2,#3,{\it #4},#5.}
\def\refinbook#1;#2;#3;#4;#5;#6.{\item{[#1]} #2, #3, #4, {\it #5},#6.} 
\def\refbook#1;#2;#3;#4.{\item{[#1]} #2,{\it #3},#4.}


\def\et{\eta}

\def\bfR{{\bf R}}

\def\calX{{\cal X}}
\def\calY{{\cal Y}}

\def\({\bigl(}
\def\){\bigr)}

\def\[{\big[}
\def\]{\big]}

\def\Ex{{\rm Ex}}

\def\lestoch{\le_{\rm I}}
\def\lerisk{\le_{\rm C}}

{
\pageno=0
\nopagenumbers
\rightline{\tt luh.arxiv.tex}
\vskip1in

\centerline{\bf Martingale Couplings and Bounds on the Tails of Probability Distributions}
\vskip0.5in

\centerline{Kyle J. Luh}
\centerline{\tt Kyle\_J\_Luh@hmc.edu}
\sk

\centerline{Nicholas Pippenger}
\centerline{\tt njp@math.hmc.edu}
\sk

\centerline{Harvey Mudd College}
\centerline{1250 Dartmouth Avenue}
\centerline{Claremont, CA 91711}
\vskip0.5in

\noindent{\bf Abstract:}
Hoeffding has shown that tail bounds on the distribution for 
 sampling from a finite population with replacement
also apply to the corresponding cases of sampling without replacement.
(A special case of this result is that binomial tail bounds apply to the corresponding hypergeometric tails.)
We give a new proof of Hoeffding's result by constructing a martingale coupling between the 
sampling distributions.
This construction is given by an explicit combinatorial procedure involving balls and urns.
We then apply this construction to create martingale couplings between other pairs of sampling distributions, both without replacement and with ``surreplacement'' (that is, sampling in which not only is the sampled individual replaced, but some number of ``copies'' of that individual are added to the population).

\vfill\eject
}

\heading{1. Introduction}

In 1963, Hoeffding [H, Section 6, Theorem 4] proved the following theorem.

\label{Theorem 1.1:}
(W.~Hoeffding)
Let the population $C$ consist of $N$ values $c_1, c_2, \ldots, c_N$.
Let $X_1, X_2, \ldots, X_n$ denote a random sample without replacement from $C$ and
let $Y_1, Y_2, \ldots, Y_n$ denote a random sample with replacement from $C$.
Let $S_n = X_1 + X_2 + \cdots + X_n$ and $T_n = Y_1 + Y_2 + \cdots + Y_n$.
Then if the function $f:\bfR\to\bfR$ is convex,
$$\Ex\[f\left(S_n\right)\] \le \Ex\[f\left(T_n\right)\].$$

Our first goal in this paper is to give a new proof of Theorem 1.1.
Our proof is based on a stochastic order relation.
The most familiar stochastic order relation is that of {\it stochastic domination}, 
which we shall denote $\lestoch$.
Stochastic domination can be defined in several equivalent ways.
Let $S$ and $T$ be real-valued random variables with finite expectations.
Then $S\lestoch T$ if $S$ and $T$ satisfy either of the following two equivalent conditions.
\medskip
\item{(I-1)}
There exists an {\it increasing coupling\/} between $S$ and $T$ (that is,
there is a random variable $(\hat{S},\hat{T})$ such that $\hat{S}$ has the same distribution as $S$,
$\hat{T}$ has the same distribution as $T$, and $\hat{S}\le \hat{T}$ with probability one).

\item{(I-2)}
For any function $f:\bfR\to\bfR$, if $f$ is increasing (that is, if $x\le y$ implies $f(x)\le f(y)$),
then $\Ex\[f(S)\]\le \Ex\[f(T)\]$.
\medskip
\noindent
(See for example M\"{u}ller and Stoyan [M, Chapter 1] or Shaked and Shathikumar 
[S2, Chapter 1].)

The stochastic order relation that is of importance in our proof is that of {\it convex domination},
which we shall denote $\lerisk$.
Convex domination can also be define in several equivalent ways.
Specifically, $S\lerisk T$ if $S$ and $T$ satisfy either of the following two equivalent conditions.
\medskip
\item{(C-1)}
There exists a {\it martingale coupling\/} between $S$ and $T$ (that is,
there is a random variable $(\hat{S},\hat{T})$ such that $\hat{S}$ has the same distribution as $S$,
$\hat{T}$ has the same distribution as $T$, and $(\hat{S},\hat{T})$ is a martingale; that is
$\Ex[\hat{T} \mid \hat{S}] = \hat{S}$).

\item{(C-2)}
For any function $f:\bfR\to\bfR$, if $f$ is convex, $\Ex\[f(S)\]\le \Ex\[f(T)\]$.
\medskip
\noindent
(See for example M\"{u}ller and Stoyan [M, Chapter 1] or Shaked and Shathikumar 
[S2, Chapter 3].)

For our proof, we shall only need the implication (C-1)$\Rightarrow$(C-2), which is easily proved as follows.
If $R$ is a random variable, we shall write $F_R(r) = \Pr[R\le r]$ for the distribution function of $R$,
so that $\Ex[R] = \int r\,dF_R(r)$.
We use the tower formula 
$\Ex[R] = \Ex\[\Ex[R\mid S]\]$ for conditional expectations, Jensen's inequality
$f\(\Ex[R]\)\le \Ex\big[f(R)\big]$ for convex $f$, and the fact that $(\hat{S},\hat{T})$ is a martingale:
$$\eqalign{
\Ex\[f(T)\]
& = \Ex\[f(\hat{T})\] \cr
&= \Ex\[\Ex\[f(\hat{T})\mid \hat{S}\]\] \cr
&= \int \int f(t)\,dF_{\hat{T}\mid\hat{S}=s}(t) \,dF_{\hat{S}}(s) \cr
&\ge \int f\left(\int t\,dF_{\hat{T}\mid \hat{S}=s}(t)\right)\,dF_{\hat{S}}(s) \cr
&= \int f\left(\Ex[\hat{T}\mid \hat{S}=s]\right)\,dF_{\hat{S}}(s) \cr
&= \int f\left(s\right)\,dF_{\hat{S}}(s) \cr
&= \Ex\[f(\hat{S})\] \cr
&= \Ex\[f(S)\]. \cr
}$$

The implication (C-1)$\Rightarrow$(C-2) shows that Theorem 1.1 is a consequence of the following proposition, which will be proved in Section 2.

\label{Proposition 1.1:}
Let the population $C$ consist of $N$ values $c_1, c_2, \ldots, c_N$.
Let $X_1, X_2, \ldots, X_n$ denote a random sample without replacement from $C$ and
let $Y_1, Y_2, \ldots, Y_n$ denote a random sample with replacement from $C$.
Let $S_n = X_1 + X_2 + \cdots + X_n$ and $T_n = Y_1 + Y_2 + \cdots + Y_n$.
Then there is a martingale coupling between $S_n$ and $T_n$.

Hoeffding used Theorem 1.1 to transfer bounds he had obtained for the tails of the distributions of sums of the independent random variables $Y_i$ to the corresponding tails for the dependent random variables $X_i$.
This transfer is possible because tail bounds typically employ a convex function, such as a quadratic or exponential, to weight large deviations from the mean more heavily than small ones.

We shall illustrate this transfer of bounds by showing how a bound on the tail of a binomially distributed random variable transfers to that of a hypergeometrically distributed random variable.
In this case, we take $N = a+b$, $c_1 = c_2 = \cdots = c_a = 1$ and $c_{a+1} = c_{a+2} = \cdots =c_{a+b} = 0$ (modeling an urn containing $a$ red balls and $b$ blue balls).
Then $S_n$  is hypergeometrically distributed (the number of red balls drawn in $n$ draws without replacement), while $T_n$ is binomially distributed (the number of red balls drawn in $n$ draws with replacement, or the number of successes in $n$ independent trials, each of which succeeds with probability $p = a/(a+b)$).

For the bound on the tail of the distribution of $T_n$ we shall use the well known method due to Chernoff [C1].
If $R$ is a random variable, we shall denote by $M_R(u) = \Ex[e^{uR}] = \int e^{ur}\,dF_R(r)$
the moment generating function of $R$.
Chernoff's bound on the upper tail of $R$ is
$$\eqalign{
\Pr[R\ge w]
&= \int_{r\ge w} dF_R(r) \cr
&\le e^{-uw}\int_{r\ge w} e^{ur}\,dF_R(r) \cr
&\le e^{-uw}\int e^{ur}\,dF_R(r) \cr
&= e^{-uw}\,M_R(u). \cr
}$$
Let $T_n$ be binomially distributed, as the number of red balls drawn in $n$ draws with replacement from an urn containing $a$ red balls and $b$ blue balls, or the number of successes in $n$ independent trials, each of which succeeds with probability $p = a/(a+b)$.
Since $M_{T_n}(u) = (pe^u + 1-p)^n$, Chernoff's bound yields
$\Pr[T_n \ge (p+q)n] \le e^{-u(p+q)n} \, (pe^u + 1-p)^n$, and minimizing this bound over $u$ yields
$$\Pr[T_n \ge (p+q)n]
\le \left(\left({p\over p+q}\right)^{p+q} \left({1-p\over 1-p-q}\right)^{1-p-q}\right)^n. \eqno(1.1)$$

We shall transfer the bound (1.1) to the corresponding tail of the corresponding hypergeometric distribution.
Let $S_n$  be hypergeometrically distributed, as the number of red balls drawn in $n$ draws without replacement from an urn containing $a$ red balls and $b$ blue balls.
The Chernoff bound on the upper tail of $S_n$ is hard to evaluate exactly (because $M_{S_n}(u)$ is a hypergeometric function, from which the distribution gets its name).
But Theorem 1.1, with the convex function $f(v) = e^{uv}$, tells us that
$$M_{S_n}(u) = \Ex[e^{uS_n}] \le \Ex[e^{uT_n}] = M_{T_n}(u).$$
Thus the Chernoff bound for $T_n$ applies to $S_n$ as well, yielding
$$\Pr[S_n \ge (p+q)n]
\le \left(\left({p\over p+q}\right)^{p+q} \left({1-p\over 1-p-q}\right)^{1-p-q}\right)^n. \eqno(1.2)$$
(Chv\'{a}tal [C2] has given a proof of the bound (1.2) by direct manipulation of sums of binomial coefficients.)

In Section 2, we shall give our construction of the martingale coupling for the proof of 
Proposition 1.1.
In Section 3, we shall apply our method to construct martingale couplings between other
pairs of distributions arising from various instances of sampling from finite populations,
without replacement, with replacement, and with ``surreplacement'' (that is, with the sampled value being replaced, together with one or more additional copies of that value).
The results of this paper first appeared in the first author's bachelor's thesis [L].
\sk

\heading{2. Proof of Proposition 1.1}

We begin with two urns.
The first urn, $\calX$, contains $N$ balls, $x_1, \ldots, x_N$.
Each of these balls is labeled with its number; that is,
ball $x_i$ is labelled $i$.
Balls will be drawn from urn $\calX$ without replacement.
The second urn, $\calY$, contains $N$ balls, $y_1, \ldots, y_N$.
Each of these balls is initially unlabeled
but will eventually be assigned a label.
Balls will be drawn from urn $\calY$ with replacement.

We now perform an infinite sequence of steps as follows.
In the course of these steps we shall define a bijective map
$\xi:\{1,\ldots,N\}\to\{1,\ldots,N\}$ and a surjective map 
$\et:\{1,2,\ldots\}\to \{1,\ldots,N\}$.
At each step, we draw a ball from urn $\calY$.
If the ball drawn is still unlabeled, we draw a ball from urn $\calX$,
we assign the label of the ball drawn from urn $\calX$ to the ball drawn from 
urn $\calY$, then replace the ball drawn from urn $\calY$ in urn $\calY$.
If the ball drawn from urn $\calY$ has already been assigned a label,
we simply replace it in urn $\calY$.
Since, with probability one, every ball in urn $\calY$ will eventually be drawn,
every ball in urn $\calY$ will eventually be assigned a label.

We define $\xi(i)$ to be the label of the $i$-th ball drawn from urn $\calX$.
Since every ball from $\calX$ is eventually drawn, and balls are drawn from $\calX$ without replacement, $\xi$ is a permutation of $\{1,\ldots, N\}$.
We define $\et(i)$ to be the label assigned to the ball drawn from urn $\calY$ at the $i$-th step
(either during the $i$-th step or at some previous step).
Since each of the labels $1,\ldots,N$ is eventually assigned to one of the balls in urn $\calY$,
$\et$  maps $\{1,2,\ldots\}$ onto $\{1,\ldots, N\}$. 

The process just described creates a coupling between $\xi$, which is uniformly distributed over all permutations of $\{1,\ldots, N\}$, and $\et$, which is a sequence $\et(1), \et(2),\ldots$ of independent random variables, each uniformly distributed over  $\{1,\ldots, N\}$.

Let $c_1, \ldots, c_n$ be real numbers.
We shall define the random variables $X_1,\ldots, X_N$ by $X_i = c_{\xi(i)}$ for $1\le i\le N$, and the random variables $Y_1, Y_2,\ldots$ by $Y_i = c_{\et(i)}$ for $i\ge 1$.
This definition creates a coupling between the sequence $X_1,\ldots, X_N$, which is distributed
as a random sample without replacement from the population $c_1, \ldots, c_n$, and
the sequence $Y_1, Y_2,\ldots$, which is distributed as a sequence of independent random samples with replacement from the same population.

Let $n$ be an integer in the range $1\le n\le N$.
We define $S_n = X_1 + \cdots + X_n$ and $T_n = Y_1 + \cdots + Y_n$.
This definition creates a coupling between $S_n$ which is distributed as the sum of a random sample of size $n$ without replacement from the population $c_1, \ldots, c_n$, and $T_n$,
which is distributed as the sum of a random sample of size $n$ with replacement from the same population.

We shall now show that $(S_n,T_n)$ is a martingale; that is, that
$$\Ex[T_n\mid S_n] = S_n. \eqno(2.1)$$
If $S_n = s$, then $c_{\xi(1)} + \cdots + c_{\xi(n)} = s$, and $\xi(1), \ldots, \xi(n)$ is equally likely to be any of the sequences satisfying this constraint.
Since any permutation of such a sequence is again such a sequence, 
we have
$$\Ex[c_{\xi(i)}\mid S_n=s] = s/n \eqno(2.2)$$
for $1\le i\le n$.
Now
$$\eqalignno{
\Ex[T_n\mid S_n=s] 
&= \Ex[Y_1\mid S_n = s] + \cdots + \Ex[Y_n\mid S_n = s] \cr
&= \Ex[c_{\et(1)}\mid S_n=s] + \cdots + \Ex[c_{\et(n)}\mid S_n=s]. &(2.3) \cr
}$$
Since each $\et(i)$ for $1\le i\le n$ is equal to one of the $\xi(1),\ldots,\xi(n)$,
each of the $n$ terms in (2.3) is equal by (2.2) to $s/n$, and thus
$\Ex[T_n\mid S_n=s] = s$.
This completes the proof of (2.1), and shows that the coupling $(S_n,T_n)$ is a martingale.
\sk

\heading{3.  Other Martingale Couplings}

In this section we shall construct martingale couplings for other pairs of probability distributions. 
(For these pairs, neither distribution has a simple moment generating function, so they do not facilitate the transfer of tail bounds in the same way as Proposition 1.1.)
The first of these pairs compares samples without replacement from two populations, one of which is a ``$k$-fold multiplication''
of the other (that is, contains $k$ ``copies'' of each individual from the other population). 

\label{Proposition 3.1:}
Let the population $C$ consist of $N$ values $c_1, c_2, \ldots, c_N$.
Let the population $D = kC$ consist of $kN$ values $d_{1,1} = \cdots = d_{1,k} = c_1, \ldots, 
d_{N,1} = \cdots = d_{N,k} = c_N$.
Let $X_1, X_2, \ldots, X_n$ denote a random sample without replacement from $C$ and
let $Y_1, Y_2, \ldots, Y_n$ denote a random sample without replacement from $D$.
Let $S_n = X_1 + X_2 + \cdots + X_n$ and $T_n = Y_1 + Y_2 + \cdots + Y_n$.
Then there is a martingale coupling between $S_n$ and $T_n$.

\label{Proof:}
We begin with two urns.
The first urn, $\calX$, contains $N$ balls, $x_1, \ldots, x_n$.
Each of these balls is labeled with its number; that is,
ball $x_i$ is labelled $i$.
Balls will be drawn from urn $\calX$ without replacement.
The second urn, $\calY$, contains $kN$ balls, $y_{1,1},\ldots,y_{1,k}, \ldots, y_{N,1},\ldots,y_{N,k}$.
Each of these balls is initially unlabeled
but will eventually be assigned a label.
Balls will be drawn from urn $\calY$ without replacement.
For $1\le i\le N$, the balls $y_{m,1},\ldots,y_{m,k}$ will be said to comprise the {\it $m$-th cohort}.

We now perform a sequence of $kN$ steps as follows.
In the course of these steps we shall define a bijective map
$\xi:\{1,\ldots,N\}\to\{1,\ldots,N\}$ and a surjective map 
$\et:\{1,\ldots,kN\}\to \{1,\ldots,N\}$.
At each step, we draw a ball from urn $\calY$.
If the ball drawn is still unlabeled, we draw a ball from urn $\calX$,
we assign the label of the ball drawn from urn $\calX$ to the ball drawn from 
urn $\calY$ and to the $k-1$ other balls in its cohort.
The ball drawn from urn $\calY$ is not replaced, and the other balls in its cohort remain in the urn.
If the ball drawn from urn $\calY$ has already been assigned a label, we
proceed to the next step.

We define $\xi(i)$ to be the label of the $i$-th ball drawn from urn $\calX$.
Since every ball from $\calX$ is eventually drawn, and balls are drawn from $\calX$ without replacement, $\xi$ is a permutation of $\{1,\ldots, N\}$.
We define $\et(i)$ to be the label assigned to the ball drawn from urn $\calY$ at the $i$-th step
(either during the $i$-th step or at some previous step).
Since each of the labels $1,\ldots,N$ is eventually assigned to one of the balls in urn $\calY$,
$\et$  maps $\{1,\ldots,kN\}$ onto $\{1,\ldots, N\}$. 

The process just described creates a coupling between $\xi$, which is uniformly distributed over all permutations of $\{1,\ldots, N\}$, and $\et$, which is uniformly distributed over maps $\et:
\{1,\ldots,kN\}$ such that $\et(h) = j$ for exactly $k$ values of $h$, for all $1\le j\le N$.

We shall define the random variables $X_1,\ldots, X_N$ by $X_i = c_{\xi(i)}$ for $1\le i\le N$, and the random variables $Y_1,\ldots, Y_{kN}$ by $Y_i = c_{\et(i)}$ for $1\le i\le kN$.
This definition creates a coupling between the sequence $X_1,\ldots, X_N$, which is distributed
as a random sample without replacement from the population $c_1, \ldots, c_N$, and
the sequence $Y_1, \ldots, Y_{kN}$, which is distributed as a sequence of independent random samples without replacement from the population $D$.
The proof this coupling is a martingale is exactly as in
the proof of Proposition 1.1.
\QED

An obvious question left open by Proposition 3.1 is whether there is a martingale coupling between sampling without replacement from population $kC$ and sampling without replacement from population $k'C$ (where $k'>k>1$, with $k$ not dividing $k'$).

Our final theorem concerns sampling with ``surreplacement'', in which 
not only is each individual drawn from a population replaced, but some number of ``copies'' of that individual are added to the population.

\label{Proposition 3.2:}
Let the population $C$ consist of $N$ values $c_1, c_2, \ldots, c_N$.
Let $X_1, X_2, \ldots, X_n$ denote a random sample without replacement from $C$ and
let $Y_1, Y_2, \ldots, Y_n$ denote a random sample with surreplacement from $C$,
whereby each individual drawn is replaced by a total of $d\ge 1$ copies of that individual.
Let $S_n = X_1 + X_2 + \cdots + X_n$ and $T_n = Y_1 + Y_2 + \cdots + Y_n$.
Then there is a martingale coupling between $S_n$ and $T_n$.
(The case $d=1$ is simply the case of sampling with replacement, dealt with in Proposition 1.1.)

\label{Proof:}
We begin with two urns.
The first urn, $\calX$, contains $N$ balls, $x_1, \ldots, x_N$.
Each of these balls is labeled with its number; that is,
ball $x_i$ is labelled $i$.
Balls will be drawn from urn $\calX$ without replacement.
The second urn, $\calY$, contains $N$ balls.
Each of these balls is initially unlabeled
but will eventually be assigned a label.
Balls will be drawn from urn $\calY$ with surreplacement.

We now perform an infinite sequence of steps as follows.
In the course of these steps we shall define a bijective map
$\xi:\{1,\ldots,N\}\to\{1,\ldots,N\}$ and a surjective map 
$\et:\{1,2,\ldots\}\to \{1,\ldots,N\}$.
At each step, we draw a ball from urn $\calY$.
If the ball drawn is still unlabeled, we draw a ball from urn $\calX$,
we assign the label of the ball drawn from urn $\calX$ to the ball drawn from 
urn $\calY$, and to $d-1$ new balls, then replace these $d$ balls in urn $\calY$.
If the ball drawn from urn $\calY$ has already been assigned a label,
we assign that label to $d-1$ new balls, then
replace these $d$ balls in urn $\calY$.
Let us consider a ball initially in urn $\calY$.
The probability that it is not drawn in the first step is $1 - 1/N$,
the probability that it is not drawn on the second step is $1 - 1/(N+(d-1))$,
and so forth, with the probability that it is not drawn on the $i$-th step being $1 - 1/(N+(i-1)(d-1)$.
Since the sum $\sum_{i\ge 1} 1/(N+(i-1)(d-1)$ diverges to infinity, the product
$\prod_{i\ge 1} \(1 - 1/(N+(i-1)(d-1)\)$ diverges to zero.
Thus, with probability one, every ball initially in urn $\calY$ will eventually be drawn,
so every ball initially in urn $\calY$ will eventually be assigned a label.
Of course, the balls added to $\calY$ are assigned labels at the times they are added.

We define $\xi(i)$ to be the label of the $i$-th ball drawn from urn $\calX$.
Since every ball from $\calX$ is eventually drawn, and balls are drawn from $\calX$ without replacement, $\xi$ is a permutation of $\{1,\ldots, N\}$.
We define $\et(i)$ to be the label assigned to the ball drawn from urn $\calY$ at the $i$-th step
(either during the $i$-th step or at some previous step).
Since each of the labels $1,\ldots,N$ is eventually assigned to one of the balls in urn $\calY$,
$\et$  maps $\{1,2,\ldots\}$ onto $\{1,\ldots, N\}$. 

The process just described creates a coupling between $\xi$, which is uniformly distributed over all permutations of $\{1,\ldots, N\}$, and $\et$, which is an sequence $\et(1), \et(2),\ldots$ of random variables, each distributed over  $\{1,\ldots, N\}$ in the way appropriate to surreplacement.
Specifically, for $i\ge 1$, the conditional probability that $\et(i)=j$, given that $\et(h)=j$ for exactly 
$k$ values of $h<i$ is $\(1+k(d-1)\)\big/\(N + (i-1)(d-1)\)$.

Let $c_1, \ldots, c_n$ be real numbers.
We shall define the random variables $X_1,\ldots, X_N$ by $X_i = c_{\xi(i)}$ for $1\le i\le N$, and the random variables $Y_1, Y_2,\ldots$ by $Y_i = c_{\et(i)}$ for $i\ge 1$.
This definition creates a coupling between the sequence $X_1,\ldots, X_N$, which is distributed
as a random sample without replacement from the population $c_1, \ldots, c_N$, and
the sequence $Y_1, Y_2,\ldots$, which is distributed as a sequence of independent random samples with surreplacement from the same population.

Let $n$ be an integer in the range $1\le n\le N$.
We define $S_n = X_1 + \cdots + X_n$ and $T_n = Y_1 + \cdots + Y_n$.
This definition creates a coupling between $S_n$ which is distributed as the sum of a random sample of size $n$ without replacement from the population $c_1, \ldots, c_n$, and $T_n$,
which is distributed as the sum of a random sample of size $n$ with surreplacement from the same population.
The proof this coupling is a martingale is exactly as in
the proof of Proposition 1.1.
\QED

An obvious question left open by Proposition 3.2 is whether there is a martingale coupling between sampling with surreplacement of $d$ copies and sampling with surreplacement of $d'$ copies from the same population, where $d'>d>1$.
\sk

\heading{4. Acknowledgment}

The research reported here was supported in part
by Grant CCF  0917026 from the National Science Foundation.
\vfill\eject

\heading{5. References}

\ref C1; H. Chernoff;
``A Measure of the Asymptotic Efficiency for Tests of a Hypothesis
Based on the Sum of Observations'';
Ann.\ Math.\ Statis.; 23 (1952) 493--507.

\ref C2; V. Chv\'{a}tal;
The Tail of the Hypergeometric Distribution;
Discr.\ Math.; 25:3 (1979) 285--287.

\ref H; W. Hoeffding;
``Probabilty Inequalities for Sums of Bounded Random Variables'';
J. Amer.\ Stat.\ Assoc.; 58:301 (1963) 13--30.

\refbook L; K. Luh;
Martingale Couplings and Bounds on the Tails of Probability Distributions;
B.~S. Thesis, Department of Mathematics, Harvey Mudd College,
May, 2011.

\refbook M; A. M\"{u}ller and D. Stoyan;
Comparison Methods for Stochastic Models and Risks;
John Wiley \& Sons, 2002.

\refbook S; M. Shaked and J. G. Shanthikumar;
Stochastic Orders;
Springer-Verlag, New York, 2007.

\bye